\documentclass{amsart}
\usepackage{amssymb}
\usepackage{graphics}
\usepackage{amsmath}

\theoremstyle{plain}
\newtheorem{lemma}{Lemma}[section]
\newtheorem{proposition}[lemma]{Proposition}
\newtheorem{theorem}[lemma]{Theorem}

\newtheorem{corollary}[lemma]{Corollary}

\theoremstyle{definition}

\newtheorem{remark}[lemma]{Remark}
\newtheorem{question}[lemma]{Question}

\newcommand\Z{{\mathbb Z}}
\newcommand\F{{\mathbb F}}
\newcommand\Q{{\mathbb Q}}
\newcommand\SL{{\operatorname{SL}}}

\begin{document}
\title[Maximal unramified $3$-extensions]{Maximal unramified $3$-extensions\\
  of imaginary quadratic fields  and $\SL_2(\Z_3)$}
\author{L. Bartholdi}
\address{\'Ecole Polytechnique F\'ed\'erale de Lausanne (EPFL),
  Institut de Math\'ematiques B (IMB),
  CH-1015 Lausanne, Switzerland}
\email{laurent.bartholdi@gmail.com}
\author{M. R. Bush}
\address{Dept.\ of Mathematics \& Statistics,
  University of Massachusetts, Amherst, MA 01003, USA}
\email{bush@math.umass.edu}
\date{February 16, 2006}
\thanks{The second author completed some of this work while he
was a graduate student at the University of Illinois (Urbana-Champaign). 
Both authors would like to thank Nigel Boston for his many insights, suggestions and support. }
\begin{abstract} 
The structure of the Galois group of the maximal unramified $p$-extension
of an imaginary quadratic field is restricted in various ways. In this paper
we construct a family of finite $3$-groups satisfying these restrictions. We
prove several results about this family and characterize them as finite
extensions of certain quotients of a Sylow pro-$3$ subgroup of $\SL_2(\Z_3)$.
We verify that the first group in the family does indeed arise as such a 
Galois group and provide a small amount of evidence that this may hold
for the other members. If this were the case then 
it would imply that there is no upper bound on the possible lengths of a finite $p$-class
tower.
\end{abstract}
\maketitle

\section{Maximal unramified $p$-extensions and Schur-$\sigma$ groups}
Let $k$ be an imaginary quadratic number field and 
$p$ be a prime.  The $p$-class tower of $k$ is the sequence of fields
\[ k = k_1 \subseteq k_2 \subseteq k_3 \subseteq \ldots \] where
$k_{n+1}$ is the maximal unramified abelian $p$-extension of $k_n$.
By Galois theory the fields $k_n$ correspond to the subgroups in the
derived series of $G = \text{Gal}\,(k^{nr,p}/k)$ where $k^{nr,p} =
\bigcup_{n \geq 1} k_n$ is the maximal unramified
$p$-extension of $k$.  If we let $Cl_p(F)$ denote the $p$-Sylow
subgroup of the class group of a number field $F$ then by class field
theory $\text{Gal}\,(k_{n+1}/k_n) \cong Cl_p(k_n)$ for $n \geq 1$. In
particular $G/[G,G] \cong \text{Gal}\,(k_{2}/k_1) \cong Cl_p(k)$ and
so is finite.


Now assume also that  $p \neq 2$.
In~\cite{KV} the notion of a Schur-$\sigma$ group is introduced.
It encapsulates various properties that the Galois group
$G$ is known to satisfy in this case.  These are:
\begin{itemize}
\item[1.] The generator rank and relation rank of $G$ (as a pro-$p$ group) are equal;
\item[2.] $G/[G,G]$ is finite;
\item[3.] There exists an automorphism $\sigma$ of order $2$ on $G$
which induces the inverse automorphism $a \mapsto a^{-1}$ on  $G/[G,G]$. 
\end{itemize}
Several structural results are proved there 
about the presentations of such groups. One consequence of their
work is that the extension $k^{nr,p}/k$ is finite only when the 
generator rank $d(G) \leq 2$. In particular all such extensions which are finite
and non-abelian  must have  $d(G)  =  2$.

In general it is exceedingly difficult to compute the Galois group
$G$.  For those examples in which the group is known to be finite the
length of the derived series is usually small. Indeed to date the
largest length observed is~$3$, see \cite{MB1}. In the
next section we will define a family of finite Schur-$\sigma$ groups
and then show that the derived length for groups in this family is
unbounded. In the last section we show that the first group in the
family is isomorphic to $\text{Gal}\,(k^{nr,p}/k)$ for several
different choices of $k$.

\section{A family of Schur $\sigma$-groups of unbounded derived length}

Let $F$ be the free pro-$3$ group on two generators $x$ and $y$.
Let $G_n$  be the Schur-$\sigma$ group defined by the pro-$3$ presentation
\[ G_n = \langle x, y \,\mid\, r_n^{-1}\sigma(r_n),\,
t^{-1}\sigma(t) \rangle \] where $r_n = x^3 y^{-3^n}$, $t = y x y
x^{-1} y$ and $\sigma : F \rightarrow F$ is the automorphism defined
by $x \mapsto x^{-1}$ and $y \mapsto y^{-1}$.  We will prove the
following result.
\begin{theorem}\label{thm:1}
For $n \geq 1$ the following hold:
\begin{itemize}
\item[(i)] $G_n$ is a finite $3$-group of order $3^{3 n + 2}$;
\item[(ii)] $G_n$ is nilpotent of class $2 n + 1$;
\item[(iii)] $G_n$ has derived length $\lfloor\log_2(3 n + 3)\rfloor$.
\end{itemize}
\end{theorem}

The remainder of this section is devoted to the proof. We first define
some auxiliary groups which are easier to study than $G_n$. Let $H_n$
be given by the pro-$3$ presentation
\[ H_n = \langle x, y \,\mid\, x^3,\,y^{3^n},\,t^{-1}\sigma(t) \rangle,\]
and let $H$ be given by the pro-$3$ presentation
\[ H = \langle x, y \,\mid\, x^3,\,t^{-1}\sigma(t) \rangle.\]
\begin{lemma}\label{lem:1}
  $G_n$ is a central cyclic extension of $H_n$; and all $H_n$'s are
  natural quotients of $H$.
\end{lemma}
\begin{proof}
  The first relation of $G_n$ is $x^6=y^{2\cdot3^n}$, so $x^6$ is
  central in $G_n$. Now the relator $x^6$ is equivalent to the
  relator $x^3$ in a 3-group; and the same argument holds for the
  relator $y^{3^n}$. It follows that the kernel of the natural map
  $G_n\to H_n$ is generated by $x^3$.
  The second assertion of the lemma is obvious.
\end{proof}

The next lemma exhibits an explicit, matrix representation of $H$. Let
$\alpha\in\Z_3$ satisfy $\alpha^2=-2$.
\begin{lemma}
  The map $\rho:H\to \SL_2(\Z_3)$, given by
  \[x\mapsto\begin{pmatrix}0 & -1\\1 & -1\end{pmatrix},\qquad
  y\mapsto\alpha\begin{pmatrix}0 & 1/2\\1 & -1\end{pmatrix},\] is an
  isomorphism between $H$ and a pro-$3$ Sylow subgroup of $\SL_2(\Z_3)$.
\end{lemma}
We recall the recursive definition of the lower $p$-central series of
a pro-$p$ group $G$: a series of closed subgroups of $G$
\[G = P_1(G) \geq P_2(G) \geq P_3(G) \geq \ldots\] defined by $P_k(G)
= P_{k-1}(G)^p [G,P_{k-1}(G)]$ for each $k \geq 1$.  Here the group on
the right-hand side is the closed subgroup generated by all $p$-th
powers of elements in $P_{k-1}(G)$, and commutators of elements from
$G$ and $P_{k-1}(G)$.

\begin{proof}
  We first claim that $\rho$ is a homomorphism. Let $\sigma':\SL_2\to
  \SL_2$ be conjugation by $(\begin{smallmatrix}-1 & 1\\0 &
    1\end{smallmatrix})$. It is then easy to check
  $\sigma'\rho=\rho\sigma$ and $\rho(x)^3=\rho(t^{-1}\sigma(t))=1$~

  We will now show simultaneously that $\rho$ is injective, and that
  its image $P$ is a pro-$3$ Sylow subgroup of $\SL_2(\Z_3)$.

  Consider the subgroup $K$ of index $3$ in $H$ that is the normal
  closure of $yx^{-1}$. It is generated by $z_i=x^{-i}yx^{i-1}$, for
  $i\in\{0,1,2\}$. Its presentation, obtained by rewriting the
  presentation of $H$ with respect to the Schreier transversal $\{e,
  x, x^2 \}$, is given by
  \[K=\langle z_0,z_1,z_2\,|\,z_iz_{i+1}^2z_i^2z_{i+1}\text{ for
  }i=0,1,2\rangle.
  \]
  The relators of $K$ may be written as
  $[z_i,z_{i+1}^{-1}]z_{i+1}^3z_i^3$; therefore, inductively, every
  element of $[K,K]$ may be written as a cube ($K$ is said to be
  ``powerful'', see~\cite{dusaut,LM1,LM2}). It follows that the lower central
  series $(\gamma_k(K))$ coincides with $(P_k(K))$ with $p=3$.  Then
  $\gamma_k(K)$ is generated, modulo $\gamma_{k+1}(K)$, by
  $\{z_i^{3^{k-1}}\}_{0\le i\le2}$.  We conclude that
  $\gamma_k(K)/\gamma_{k+1}(K)$ has rank at most $3$.

  Recall that $\SL_2(\Z_3)$ has congruence kernels
  $N_k=1+3^kM_2(\Z_3)$. The lower central series of $N_1$ is given by
  $\gamma_{k}(N_1)=N_k$, and the rank of $N_k/N_{k+1}$ is $3$.

  All the claims will follow if we show that $\{\rho(z_i^{3^{k-1}})\}$
  spans $N_k/N_{k+1}$ for all $k\ge1$; indeed then $\rho(K)=N_1$, and
  since the ranks along the lower central series of $K$ are bounded by
  $3$, they must equal $3$ and $\rho$ is then injective. We compute:
  \[\rho(z_0^{3^{k-1}})=\begin{pmatrix}\alpha^{-3^{k-1}} & 0\\ 0 &
    \alpha^{3^{k-1}}\end{pmatrix};\] and
  $\alpha^{3^{k-1}}\in1+3^k\Z_3\setminus1+3^{k+1}\Z_3$; or, in
  other words, $\alpha$ is a topological generator of the
  torsion-free part of $\Z_3^\times$.
  Similarly,
  \begin{gather*}
    \rho(z_1^{3^{k-1}})=\begin{pmatrix}\alpha^{3^{k-1}} &
      \alpha^{-3^{k-1}}-\alpha^{3^{k-1}}\\ 0 &
      \alpha^{-3^{k-1}}\end{pmatrix},\\
    \rho(z_2^{3^{k-1}})=\begin{pmatrix}\alpha^{3^{k-1}} & 0\\
      \alpha^{-3^{k-1}}-\alpha^{3^{k-1}} &
      \alpha^{-3^{k-1}}\end{pmatrix},
  \end{gather*}
  and the off-diagonal entries are in $3^k\Z_3\setminus3^{k+1}\Z_3$.

  We conclude by considering $P=N_1\langle \rho(x)\rangle$, which is a
  pro-$3$ Sylow subgroup of $\SL_2(\Z_3)$, and noting that $\rho(H)=P$ since
  they have isomorphic index-$3$ subgroups $K$ and $N_1$
  respectively.
\end{proof}

\begin{remark}
  (i) The proof is similar to a construction of presentations of
  congruence kernels in~\cite{AndozCvetk}.

  (ii) The following simple and general argument was generously
  communicated to us by Nigel Boston and Jordan Ellenberg, see 
  \cite{BostEll}. Suppose $f:
  T\to U$ is a surjective homomorphism of pro-$p$ groups such that
  $H^1f: H^1(U,\F_p)\to H^1(T,\F_p)$ is an isomorphism $H^2f:
  H^2(U,\F_p)\to H^2(T,\F_p)$ is injective. Then $f$ is an
  isomorphism.

  We may apply it to $T=K$ and $U=N_1$. It is not difficult to show
  that $f$ is surjective and that $H^1f$ is an isomorphism. Now the
  cup product map $\bigwedge^2 H^1(U,\F_p)\to H^2(U,\F_p)$ is an
  isomorphism, because $U$ is uniform; so to prove injectivity of
  $H^2f$ it suffices to show that $\bigwedge^2 H^1(T,\F_p)\to
  H^2(T,\F_p)$ is injective; this holds because $T/\Phi(\Phi(T))$ is
  abelian.
\end{remark}

We may now identify $H_n$ with an appropriate quotient of $P$:
\begin{lemma}
  $H_n$ is the quotient of $P$ by the subgroup of matrices
  $(\begin{smallmatrix} a & b\\c & d\end{smallmatrix})$ satisfying the
  congruences
  \begin{gather*}
    a,d\equiv 1\pmod{3^n}\\
    b,c\equiv 0\pmod{3^n}\\
    a+d\equiv 2\pmod{3^{n+2}}\\
    a+b\equiv 1\pmod{3^{n+1}}\\
    a+b-c\equiv 1\pmod{3^{n+2}}
  \end{gather*}
\end{lemma}
\begin{proof}
  This amounts to computing the normal closure $R$ of $\rho(y^{3^n})$ in
  $P$.  We have $y^3$ conjugate to $z_2z_1z_0$, which implies $\rho(y^{3^n})\equiv
  \rho(z_0^{3^{n-1}}z_1^{3^{n-1}}z_2^{3^{n-1}})$ in
  $N_{n}/N_{n+1}$, and so the intersection of $R$ with $N_{n}/N_{n+1}$
  is one-dimensional.

  Then, taking commutators with $z_i$, we have $[y^3,z_i]\equiv
  z_{i-1}^{3^n}z_i^{3^n}/z_i^{3^n}z_{i+1}^{3^n}\equiv
  z_{i-1}^{3^n}/z_i^{3^n}$ in $N_{n+1}/N_{n+2}$; so the intersection
  of $R$ with $N_{n+1}/N_{n+2}$ is two-dimensional.

  Finally, taking a commutator again, we have
  \[[z_{i-1}^{3^n}/z_i^{3^n},z_{i+1}]\equiv z_{i+1}^{-3^{n+1}}\] in
  $N_{n+2}/N_{n+3}$, so the intersection of $R$ with $N_{n+2}/N_{n+3}$
  is three-dimensional, and the same holds for $N_{k}/N_{k+1}$ for all
  $k\ge n+2$.

  Writing the matrices $\rho(z_i^{3^n})$ then proves the lemma.
\end{proof}

Finally, we identify better the relation between $G_n$ and $H_n$:
\begin{lemma}\label{lem:4}
  The kernel of the natural map $G_n\to H_n$ is cyclic of order $3$.
\end{lemma}
\begin{proof}
  The kernel is cyclic by Lemma~\ref{lem:1}. The order of $y^{3^n}$ in
  $G_n$ is at most $3$, since
  $y^{3^n}=z_0^{3^{n-1}}z_1^{3^{n-1}}z_2^{3^{n-1}}$ and the relations
  in $H_n$ force cubes of $z_i^{3^{n-1}}$ to be commutators, and
  therefore to vanish in any central extension.

  On the other hand, no relation in $G_n$ prevents $y^{3^n}$ from
  being
  non-trivial in $G_n$; therefore the order of $y$ is precisely
  $3^{n+1}$, and $|G_n|=3|H_n|$.
\end{proof}

We are finally ready to prove the main theorem of this section:
\begin{proof}[Proof of Theorem~\ref{thm:1}]
  (i) We have $|G_n|=3|H_n|$ by Lemma~\ref{lem:4}, and
  $|H_n|=3^{3n+1}$ because the normal closure $R_n$ of $y^{3^n}$ in
  $H$ has index $3^{3n}$ in $N_1$, and therefore has index $3^{3n+1}$
  in $P$.

  (ii) We first compute the lower central series of $H_n$. It is
  obtained from that of $P$, as follows: $\gamma_1(P)=P$; and for
  $k\ge1$,
  $\gamma_{2k}(P)=N_{k+1}\langle(z_0/z_1)^{3^k},(z_1/z_2)^{3^k}\rangle$
  and $\gamma_{2k+1}(P)=N_{k+1}\langle(z_0z_1z_2)^{3^k}\rangle$. The
  last index $k$ such that $R_n$ is not contained in $N_k$ is $n+1$,
  so the nilpotency class of $H_n$ is $2n+1$. Finally, the action by
  conjugation of $x$ on $y^{3^n}\equiv(z_0z_1z_2)^{3^{n-1}}$ is
  trivial, so the nilpotency class of $G_n$ is the same as that of
  $H_n$, namely $2n+1$; the last quotient
  $\gamma_{2n+1}(G_n)/\gamma_{2n+2}(G_n)=\langle x,y^{3^n}\rangle$.

  (iii) The derived length of $G_n$ can also be obtained from the
  derived series of $P$: one has $P^{(2k)}=\gamma_{(2^{2k+2}-1)/3}(P)$
  and $P^{(2k+1)}=\gamma_{(2^{2k+3}-2)/3}(P)$ using
  $[N_k,N_\ell]=N_{k+\ell}$, which comes from the identity
  \[[1+3^mA,1+3^nB]\equiv1+3^{m+n}(AB-BA)\]
  and the fact that the Lie algebra $sl_2$ is simple.

  By~(ii), we have $\gamma_{2n+1}(P)>R_n>\gamma_{2n+2}(P)$, so
  $P^{(k)}>R_n>P^{(k+1)}$ for $k=\lfloor\log_2(3 n + 3)\rfloor$. The
  same argument as above shows that the derived length of $G_n$ is
  equals that of $H_n$.
\end{proof}

\begin{remark}
  The groups $G_n$ are finite $3$-groups with the same number of
  relations as generators in their pro-$3$ presentations. It is an
  open question  as to whether or not this
  implies that such groups must have an \emph{abstract} presentation with
  equal numbers of generators and relations. Finite groups with this
  latter property are said to have {\em deficiency zero}.
  It is also open whether or not there exist abstract groups of
  deficiency zero with arbitrarily large derived length.  To date the
  maximum length achieved is $6$ (see~\cite{Kenne}). If $G_n$ has
  deficiency zero then this question would be resolved.

We note that examples similar to ours have appeared in the literature
before.  In~\cite{AndozCvetk} a family of finite $3$-generated
$p$-groups (for odd prime $p$) with increasing nilpotency class and
derived length is constructed.  However our family of groups are the
first $2$-generated candidates to appear in the literature, as far as
we know.
\end{remark}

\section{Explict computations of $\text{Gal}\,(k^{nr,3}/k)$}
In~\cite{BG1} and~\cite{MB1} the $p$-group generation algorithm
is used to compute the Galois groups of several $p$-extensions
with restricted ramification. Here we use it to verify that 
$\text{Gal}\,(k^{nr,3}/k) \cong G_1$
for several different imaginary quadratic fields $k$.
For the reader's convenience we recall some definitions and 
give a brief description of the method.

Recall that $(P_k(G))$ denotes the lower $p$-central series of $G$. If
$G$ is a finite $p$-group then this series is finite and the
smallest $c$ such that $P_c(G) = \{1\}$ will be called the {\em
  $p$-class} of $G$. A $p$-group $H$ is called a {\em descendant} of
$G$ if $H / P_c(H) \cong G$ where $c$ is the $p$-class of $G$. It is
an {\em immediate descendant} if it has $p$-class $c + 1$. The
$p$-group generation algorithm~\cite{Obrien} finds representatives (up
to isomorphism) of all the immediate descendants of a given finite
$p$-group $G$. Starting with the elementary abelian $p$-group on $d$
generators (for some fixed $d$) the algorithm allows one to compute a
tree containing all finite $d$-generated $p$-groups. The $p$-class of
a group determines the level of the tree in which it occurs.

For the Galois groups we are interested in we have additional information
about the maximal abelian quotients of various subgroups of small index.
This information is obtained by computing class groups of various extensions
and applying the Artin reciprocity isomorphism from class field theory.
This information is sometimes sufficient to eliminate all but finitely
many groups from the tree of descendants described above, in which case
we are left with a finite number candidates for the Galois group.
A more precise formulation of the method and several examples
in the case $p = 2$ can be found in~\cite{MB1}.

In the case $p = 3$ we have obtained the following result using
the symbolic computation package \textsc{Magma}~\cite{magma}. (Note: 
we describe abelian groups by listing the orders of their cyclic components.
So for instance $[3,3]$ is the direct product of a cyclic group of order 3 
with itself.)
\begin{proposition}\label{prop_identify3grp}
  Let $G$ be a pro-$3$ group and suppose $G/[G,G] \cong [3,3]$, then $G$ has four
  closed subgroups of index 3.  If these four subgroups have maximal abelian
  quotients $[3,9]$, $[3,9]$, $[3,9]$ and $[3,3,3]$, then
  $G$ is a finite $3$-group.
\end{proposition}
In fact after starting the $p$-group generation algorithm on the $2$-generated
elementary abelian $3$-group $[3,3]$ with the restrictions described in the proposition
the computation terminates having found two candidates for $G$. 
These will be denoted by $Q_1$ and $Q_2$ and
 are generated by $\{x_i\}_{i = 1}^5$ subject to the
following power-commutator presentations.
\begin{align*}
  (Q_1) \qquad\qquad  &x_1^3 = x_4     &[x_2 , x_1] &= x_3 \\ 
    &x_2^3 = x_4     &[x_3 , x_1] &= x_4 \\
   &  &[x_3 , x_2] &= x_5
\end{align*}
\begin{align*}
  (Q_2) \qquad\qquad  &x_1^3 = x_4^2     &[x_2 , x_1] &= x_3 \\ 
    &    &[x_3 , x_1] &= x_4 \\
   &  &[x_3 , x_2] &= x_5
\end{align*}
$Q_1$and $Q_2$ are the groups $(243,5)$ and $(243,6)$ respectively in
\textsc{Magma}'s or \textsc{Gap}'s {\tt SmallGroups}
database~\cite{magma,gap}.
\begin{remark}
Note that in these power-commutator presentations if a power
$x_r^3$ or commutator $[x_r, x_s]$ does not occur on the left-hand
side of the given relations then it is assumed to be trivial.
\end{remark}

\begin{corollary}
All discriminants $-50000 \leq d \leq 0$ for which
the field $k = \Q(\sqrt{d})$ has $\text{Gal}\,(k^{nr,3}/k) \cong G_1$
are contained in the following list:
$d = -4027$, $-8751$, $-19651$, $-21224$, $-22711$, $-24904$, $-26139$, $-28031$, $-28759$, $-34088$, 
$-36807$, $-40299$, $-40692$, $-41015$, $-42423$, $-43192$,
$-44004$, $-45835$, $-46587$, $-48052$, $-49128$, and $-49812$.
\end{corollary}
\begin{proof}
For each of these fields $Cl_3(k) \cong [3,3]$. 
\textsc{Magma}'s class field theory package can be used to construct and verify that the
four unramified extensions $\{k_i\}_{i = 1}^4$ of degree 3 over $k$ have 
$Cl_3(k_i) \cong [3,9]$ for three choices of $i$, and $Cl_3(k_i) \cong [3,3,3]$
for the remaining choice. By Proposition~\ref{prop_identify3grp}, 
$\text{Gal}\,(k^{nr,3}/k)$ is isomorphic to $Q_1$ or $Q_2$. The group $Q_2$ has non-trivial Schur
multiplier and hence can be eliminated (see~\cite{KV}) leaving $Q_1$ 
as the only possibility. One can verify (by hand or by machine computation)
that $G_1$ satisfies the conditions in Proposition~\ref{prop_identify3grp}. Hence
we must also have $Q_1 \cong G_1$.
\end{proof}
\begin{remark}
Since $G_1$ has derived length $2$ the fields described in the corollary all have 
$3$-class towers of length $2$.  In~\cite{ST} the field $\Q(\sqrt{-4027})$ is shown to 
have $3$-class tower of length 2 by different means.
\end{remark}

The following question remains to be answered:
\begin{question}
   Is it possible, for all $n \geq 1$, to find an
  imaginary quadratic field $k$ such that $\text{Gal}\,(k^{nr,3}/k)
  \cong G_n$?
\end{question}
If the answer is yes then this would imply that the lengths of finite
$p$-class towers are unbounded. As a first step towards answering this question
one might compute the abelian quotient invariants
of the index 3 subgroups in $G_n$ for various $n \geq 2$, and then search for fields $k$
which have unramified extensions with matching $3$-class groups. Using standard
methods \cite{HJ,Johnson} one can compute the abelian quotient invariants and 
it turns out that the result is independent of $n$. More precisely one obtains 
the following proposition.
\begin{proposition}
For $n \geq 2$ the group $G_n$ has four subgroups of index $3$ with 
abelian quotient invariants $[3,9]$, $[3,3,3]$, $[3,3,3]$, and $[3,3,3]$.
\end{proposition}
Moreover when one looks for examples of imaginary quadratic fields
having unramified extensions with matching $3$-class groups they seem
relatively easy to find.  The discriminants $d$ with $|d| < 50000$
for which there is a match are $d = -3896$, $-6583$, $-23428$,
$-25447$, $-27355$, $-27991$, $-36276$, $-37219$, $-37540$, $-39819$,
$-41063$, $-43827$, $-46551$.

At this point it becomes difficult to make further progress. Clearly one cannot use the
abelian quotient invariants of the index $3$ subgroups to separate out any of the groups
$G_n$ for $n \geq 2$ as we did with $G_1$. If one restricts attention to the smallest groups
$G_2$ and $G_3$ then differences in the abelian quotient invariants only show
up when one looks at subgroups of index at least $27$. The corresponding class
group computations that one would need to carry out do not seem feasible
currently.

\end{document}